\documentclass[12pt]{article} 
\usepackage{amsmath}
\usepackage{amssymb} 
\usepackage{amscd}
\usepackage{amsthm}
\usepackage[utf8]{inputenc}
\usepackage[russian,english]{babel}
\usepackage{pb-diagram}
\usepackage{hyperref}
\def\im{{\mbox{Im}}}

\def\ker{{\mbox{Ker}}}

\def\cala{{\cal A}}

\def\cals{{\cal S}}
\def\calc{{\cal C}}

\def\bbbone{\mbox{\rm 1\hspace {-.6em} l}}

\numberwithin{equation}{section}

\begin{document}
\enlargethispage{3cm}

\thispagestyle{empty}
\begin{center}
{\bf TENSOR PRODUCT OF N-COMPLEXES}
\end{center} 
\begin{center}
{\bf AND GENERALIZATION}
\end{center}
\begin{center}
{\bf  OF GRADED DIFFERENTIAL ALGEBRAS}
\end{center}
   
\vspace{0.3cm}

\begin{center}
Michel DUBOIS-VIOLETTE
\footnote{Laboratoire de Physique Th\'eorique, UMR 8627\\
Universit\'e Paris XI,
B\^atiment 210\\ F-91 405 Orsay Cedex\\
Michel.Dubois-Violette$@$u-psud.fr} 
\end{center}
\vspace{0,1cm}

\begin{center}
\begin{otherlanguage}{russian}
 на Иван Тодоров\\
\end{otherlanguage}
{\sl Dédié à Ivan Todorov}
\end{center}
 \vspace{0,1cm}
 \begin{center}
 {\sl ``C'est une chose de toute éternité que l'amitié intellectuelle"},\\
  Auguste Villiers de l'Isle-Adam, {\sl L'Intersigne}.
\end{center}
\begin{abstract}
It is known that the notion of graded differential algebra coincides with the notion of monoid in the monoidal category of complexes. By using the monoidal structure introduced by M. Kapranov for the category of $N$-complexes we define the corresponding generalization of graded differential algebras as the monoids of this category. It turns out that this generalization coincides with the notion of graded\linebreak[4] $q$-differential algebra which has been previously introduced and studied.
\end{abstract}

\vfill

\noindent LPT-ORSAY 09/84
\newpage

\section{Introduction}
In the setting of graded algebras over $\mathbb C$, the notion of graded $q$-differential algebra has been introduced in  \cite{mdv-ker:1996}. This notion stabilizes in the more general setting of \cite{mdv:1998a}. It appears there as a generalization of the notion of graded differential algebra in which the differential is replaced by a $N$-differential, i.e. a homogeneous linear mapping $d$ of degree 1 satisfying $d^N=0$ (instead of $d^2=0$), such that one has a twisted Leibniz rule (instead of the antiderivation property).\\
Lots of examples have been introduced and analyzed in \cite{mdv-ker:1996} and 
\cite{mdv:1998a}.\\

Our aim in this article is to discuss the naturality of this notion. To understand this point, we recall that there is a class of categories called monoidal categories (or tensor categories) in which there are natural notions of algebras (or monoids). Such a monoidal category $\calc$ is equipped with a functor 
\[
\otimes : \calc \times \calc \rightarrow \calc
\]
called the tensor product satisfying associativity modulo natural isomorphisms subject to some axioms. We refer to \cite{mac:1971} for the precise definition.\\

The monoidal categories that will be considered here are categories of modules over a unital commutative ring $\mathbb K$ in which the monoidal structure, i.e. the tensor product, will always induce the ordinary tensor product on the underlying $\mathbb K$-modules. For instance one may consider the category of all $\mathbb K$-modules equipped with the usual tensor product, the category of $(\mathbb Z$-)graded $\mathbb K$-modules equipped with the usual graded tensor product or the category of (cochain) complexes of $\mathbb K$-modules equipped with the usual tensor product of complexes which is recalled in the next section.\\

A monoid (or an algebra) in such a monoidal category $\calc$ is an object $\cala$ of $\calc$ endowed with a product
\[
\mu : \cala\otimes \cala\rightarrow \cala
\]
which is associative. In the above examples, a monoid of the category of $\mathbb K$-modules is simply an associative $\mathbb K$-algebra, a monoid of the category of graded $\mathbb K$-modules is an associative graded $\mathbb K$-algebra while a monoid of the category of (cochain) complexes of $\mathbb K$-modules is a graded differential $\mathbb K$-algebra (see in the next section).\\

For the category $\calc_N$ of (cochain) $N$-complexes of $\mathbb K$-modules there is a monoidal structure (i.e. a tensor product) introduced by M. Kapranov in \cite{kap:1991}. To define it, one must assume that $\mathbb K$ has a primitive $N$-th root of the unit $q$ in a general sense (see below) and then, up to the choice of such a $q\in \mathbb K$, the monoidal structure is essentially unique. It turns out that the monoids of $\calc_N$ are precisely the graded $q$-differential algebras over $\mathbb K$ as defined in \cite{mdv:1998a}. \\

In this article we give a detailed description of these points. In the next section we recall the monoidal structure of the category of complexes and the definition of graded differential algebras as the monoids of this category. In Section 3 we introduce the category $\calc_N$ of $N$-complexes and discuss some of its properties. In Section 4 we introduce the monoidal structure of $\calc_N$ and we identify its monoids as the graded $q$-differential algebras. In conclusion we indicate that this notion is probably not the end of the story concerning the ``$N$-generalization" of graded differential algebras.\\

Throughout this article, by a complex or a $N$-complex we always mean a cochain complex or a cochain $N$-complex, that is the differential or the $N$-differential is of degree 1 (degree -1 corresponding to chain complexes or chain $N$-complexes).

\section{Preliminary : From the tensor product of complexes to graded differential algebras}

In this article, $\mathbb K$ denotes a unital commutative ring. By {\sl a complex of $\mathbb K$-modules} or simply {\sl a complex} when no confusion arises, we here mean a graded $\mathbb K$-module
\begin{equation}
C=\oplus_{n\in \mathbb Z} C^n
\label{eq1.1}
\end{equation}
equipped with an endomorphism $d$ which is homogeneous of degree 1, i.e.
\begin{equation}
d(C^n)\subset C^{n+1},\ \ \ \forall n\in \mathbb Z
\label{eq1.2}
\end{equation}
and is such that $d^2=0$, i.e.
\begin{equation}
d^2(x)=0
\label{eq1.3}
\end{equation}
for any $x\in C$.\\
There is a natural tensor product $C_0\otimes C_1$ of complexes $C_0$ and $C_1$ defined by
\begin{equation}
(C_0\otimes C_1)^n=\oplus_{r+s=n} C^r_0\otimes C^s_1
\label{eq1.4}
\end{equation}
and
\begin{equation}
d(x_0\otimes x_1)=d(x_0)\otimes x_1+(-1)^nx_0\otimes d(x_1)
\label{eq1.5}
\end{equation}
for $x_0\in C^n_0$ and $x_1\in C_1$. It is of course understood that on the right-hand side of (\ref{eq1.4}) the tensor products $C^r_0 \otimes C^s_1$ are the tensor products over $\mathbb K$ of $\mathbb K$-modules.\\
This defines a monoidal structure \cite{mac:1971} for the category of complexes which is essentially unique extending the canonical one of the category of graded modules. The graded differential algebras (over $\mathbb K$) are then defined as the monoids of this monoidal category.\\
In other words a graded differential algebra is a complex $\cala$ equipped with a product
\[
\mu :\cala\otimes \cala\rightarrow \cala
\]
which is associative. This means that $\cala$ is a graded algebra with product 
\[
(x,y)\mapsto xy=\mu(x\otimes y)
\]
which is associative of degree 0 and that, in view of (\ref{eq1.5}) one has
\begin{equation}
d(xy)=d(x)y+(-1)^n xd(y)
\label{eq1.6}
\end{equation}
for $x\in \cala^n$ and $y\in \cala$.\\
One thus recovers the usual notion of graded differential algebra. It is our aim in the following to use the monoidal structure introduced by M. Kapranov \cite{kap:1991} for the category of $N$-complexes to extract the corresponding generalization of graded differential algebras and to give some of their properties. It turns out that this generalization coincides with the notion of graded $q$-differential algebra introduced and studied in \cite{mdv-ker:1996}, 
\cite{mdv:1998a}.

\section{The category of $N$-complexes}

Throughout this article $N$ denotes an integer greater or equal to 2, i.e. $N\in \mathbb N$ with $N\geq 2$.

\subsection{$N$-complexes}
 A $N$-{\sl complex of} $\mathbb K$-modules or simply a $N$-{\sl complex} when no confusion arises is a graded $\mathbb K$-module
 \begin{equation}
 C=\oplus_{n\in \mathbb Z} C^n
 \label{eq2.1}
 \end{equation}
 equipped with an endomorphism $d$ which is homogeneous of degree 1 and satisfies $d^N=0$, i.e. one has
 \begin{equation}
 d(C^n)\subset C^{n+1},\ \ \ \forall n\in \mathbb Z
 \label{eq2.2}
 \end{equation}
 and
 \begin{equation}
 d^N(x)=0
 \label{eq2.3}
 \end{equation}
 for any $x\in C$. The endomorphism $d$ is refered to as the {\sl $N$-differential} of $C$.\\
 An {\sl homomorphism of $N$-complexes} of $C$ into $C'$ is a $\mathbb K$-linear mapping\linebreak[4] $\alpha:C\rightarrow C'$ which is homogeneous of degree 0 and such that
 \begin{equation}
 \alpha(d(x))=d(\alpha(x))
 \label{eq2.4}
 \end{equation}
for any $x\in C$, (we denote by the same symbol $d$ the $N$-differential of all $N$-complexes when no confusion arises).\\ 
The $N$-complexes and their homomorphisms form a category $\calc_N$ which is an abelian category, (see in \cite{car-eil:1973}, \cite{mac:1971} or \cite{wei:1994} for the definition of an abelian category). $N$-complexes occur in various contexts, see e.g. \cite{may:1942a}, \cite{may:1942b}, \cite{kap:1991}, \cite{mdv-ker:1996}, \cite{mdv:1996b}, \cite{kas-wam:1998}, \cite{mdv:1998a}, 
 \cite{mdv-tod:1997}, \cite{mdv-tod:1999}, \cite{mdv-hen:1999}, \cite{mdv:2002},  \cite{mdv-hen:2002}, \cite{ber-mdv-wam:2003},  \cite{ber-mar:2006},  \cite{ac-mdv:2002b},  \cite{ac-mdv:2007}, \cite{cib-sol-wis:2007} and \cite{wam:2001}.
 
\subsection{Cohomology}

Given a $N$-complex of $\mathbb K$-modules $C$, there are $N-1$ graded $\mathbb K$-modules $H_{(k)}(C)$ for $1\leq k\leq N-1$ which are the generalization of the cohomology and were called {\sl generalized cohomology of} $C$ in \cite{mdv:1998a} and {\sl amplitude cohomology of } $C$ in \cite{cib-sol-wis:2007}. They are defined by
\[
H_{(k)}(C)=\ker (d^k)/\im (d^{N-k})
\]
for $k\in \{1,\dots, N-1\}$. One has
\begin{equation}
H_{(k)}(C)=\oplus_{n\in \mathbb Z} H^n_{(k)}(C)
\label{eq2.5}
\end{equation}
where
\begin{equation}
H^n_{(k)}(C)=\{x\in C^n\vert d^k(x)=0\}/d^{N-k}(C^{n-N+k})
\label{eq2.6}
\end{equation} 
 for $n\in \mathbb Z$, $k\in \{1,\dots,N-1\}$.\\
 
 The $H_{(k)}(C)$, $k\in \{1,\dots,N-1\}$ are connected by by the exact hexagons \cite{mdv-ker:1996} \cite{mdv:1996b}, \cite{mdv:1998a}
 
 \[ \begin{diagram} \node{}
\node{H_{(\ell+m)}(C)} \arrow{e,t}{[d]^m} \node{H_{(\ell)}(C)}
\arrow{se,t}{[i]^{N-(\ell+m)}} \node{} \\ \node{H_{(m)}(C)}
\arrow{ne,t}{[i]^\ell}
\node{} \node{} \node{H_{(N-m)}(C)} \arrow{sw,b}{[d]^\ell} \\ \node{}
\node[1]{H_{(N-\ell)}(C)} \arrow{nw,b}{[d]^{N-(\ell+m)}}
\node{H_{(N-(\ell+m))}(C)}
\arrow{w,b}{[i]^m} \node{} \end{diagram} \]
for $\ell,m\in \{1,\dots,N-1\}$ such that $\ell+m\leq N-1$, where $[i]$ and $[d]$ are respectively induced by the canonical inclusions $i:\ker(d^n)\rightarrow \ker(d^{n+1})$ for $1\leq n\leq N-1$ and the mappings $x\mapsto dx$ of $\ker(d^n)$ into $\ker(d^{n-1})$ for $1\leq n\leq N-1$. Notice that $[d]$ is homogeneous of degree 1 while $[i]$ is of degree 0.
  
For any short exact sequence
 \[
 0\rightarrow C_1\stackrel{\alpha}{\rightarrow} C_2\stackrel{\beta}{\rightarrow}C_3\rightarrow 0
 \]
 of $N$-complexes there are connecting homomorphisms 
 \[
 \partial:H_{(k)}(C_3)\rightarrow H_{(N-k)}(C_1)
 \]
  for $k\in\{1,\dots,N-1\}$ such that the $H_{(n)}(C_i)$ and the $H_{(N-n)}(C_i)$ are connected by the exact hexagons \cite{mdv:1998a}
  
  \[ \begin{diagram} \node{}
\node{H_{(n)}(C_2)} \arrow{e,t}{\beta_\ast} \node{H_{(n)}(C_3)}
\arrow{se,t}{\partial} \node{} \\ \node{H_{(n)}(C_1)} \arrow{ne,t}{\alpha_\ast}
\node{} \node{} \node{H_{(N-n)}(C_1)} \arrow{sw,b}{\alpha_\ast} \\ \node{}
\node[1]{H_{(N-n)}(C_3)} \arrow{nw,b}{\partial} \node{H_{(N-n)}(C_2)}
\arrow{w,b}{\beta_\ast} \node{} \end{diagram} \]
 for $n\in \{1,\dots,N-1\}$, where $\alpha_\ast$ and $\beta_\ast$ are induced by $\alpha$ and $\beta$. Notice that $\partial$ is homogeneous of degree 1 while $\alpha_\ast$ and $\beta_\ast$ are of degree 0 (as $\alpha$ and $\beta$). 
  
 The $H_{(k)}$, $k\in \{1,\dots,N-1\}$ are covariant functors from $\calc_N$ to the category of graded $\mathbb K$-modules and the construction of the connecting homomorphism $\partial$ is also functorial as explained in \cite{mdv:1998a}.

 \section{Monoidal structure on $\calc_N$}
 
 In order to define the tensor product of $N$-complexes we need some assumptions on the ring $\mathbb K$ which we now explain.
 
 \subsection{Basic assumption on $\mathbb K$}
 
 Let $q$ be an element of the ring $\mathbb K$. One associates to $q$ a mapping\linebreak[4] $[\bullet]_q:\mathbb N\rightarrow \mathbb K,n\mapsto [n]_q$ defined by setting 
 \[
 [0]_q=0
 \]
 and
 \begin{equation}
 [n]_q=1+\dots + q^{n-1}=\sum^{n-1}_{k=0} q^k,\ \ \ \forall n\geq 1
 \label{eq3.1}
 \end{equation}
 where by convention $q^0=1$. For $n\geq 1$, one defines the $q$-factorial $[n]_q!$ by $\prod^n_{k=1} [k]_q$ and for $n$ and $m$ with $n\geq 1$ and $n\geq m\geq 0$, one defines inductively the $q$-binomial coefficients
$ \left[\begin{array}{c}n\\ m \end{array}\right]_q\in \mathbb K$ by setting
 \begin{equation}
 \left[\begin{array}{c}n\\ 0 \end{array}\right]_q=
 \left[\begin{array}{c}n\\ n \end{array}\right]_q = 1
 \label{eq3.2}
 \end{equation}
 and
  \begin{equation}
 \left[\begin{array}{c}n\\ m \end{array}\right]_q
 + q^{m+1} \left[\begin{array}{c}n\\ m+1 \end{array}\right]_q = 
 \left[\begin{array}{c}n+1\\ m+1 \end{array}\right]_q 
 \label{eq3.3}
 \end{equation}
 for $n-1\geq m\geq 0$.\\
 We now make the following basic assumption (A) on the ring $\mathbb K$\\

(A) $\left[\begin{array}{l}
\mathbb K \ \text{\sl has a distinguished element}\ q\\
\text{\sl such that}\ [N]_q=0\ \text{\sl and such that}\\
\text{\sl the}\  [n]_q\ \text{\sl are invertible in}\ \mathbb K\ \text{\sl for}\ N-1\geq n\geq 1.
\end{array}
\right.$\\

\noindent Assumption (A) has been repetedly used in \cite{kas-wam:1998} and in \cite{mdv:1998a} and means that the ring $\mathbb K$ is equipped with a distinguished primitive N$^{\text{th}}$ root of the unit.\\
 Standard examples of pairs $(\mathbb K,q)$ satisfying (A) are $(\mathbb C,e^{\frac{2\pi i}{N}})$ and, when $N$ is a prime number,  $(\mathbb Z/N\mathbb Z,1)$. Notice that for $N=2$ Assumption (A) is immaterial and $q=-1$.
 
 \subsection{Tensor product of $N$-complex}
 
 Here and in the sequel of this paper it is assumed that Assumption (A) is satisfied by the pair $(\mathbb K,q)$.\\
 Following \cite{kap:1991}, we define the tensor product $C_0\otimes C_1$ of $N$-complexes $C_0$ and $C_1$ by
 \[
 (C_0\otimes C_1)^n=\oplus_{r+s=n} C^r_0 \otimes C^s_1
 \]
 and
 \begin{equation}
  d(x_0\otimes x_1)=d(x_0)\otimes x_1+q^n x_0\otimes d(x_1)
  \label{eq3.4}
 \end{equation}
 for $x_0\in C^n_0$ and $x_1\in C_1$.\\
 The very reason why this definition works is that, by using the above definition of the $q$-binomial coefficients one obtains by induction on $n$
 \begin{equation}
 d^k(x_0\otimes x_1)=\sum^k_{p=0} q^{n(k-p)}\left[\begin{array}{c}k\\ p \end{array}\right]_q d^p (x_0) \otimes d^{k-p}(x_1)
 \label{eq3.5}
 \end{equation}
 for $x_0\in C^n_0$ and $x_1\in C_1$ and that it follows from $[N]_q=0$ and from the invertibility of the $[p]_q$ for $1\leq p\leq N-1$ that one has
 \begin{equation}
 \left[\begin{array}{c}
 N\\ p
 \end{array}\right]_q=0
 \label{eq3.6}
 \end{equation}
 for $p\in \{1,\dots,N-1\}$. Therefore one has
 \begin{equation}
 d^N(x_0\otimes x_1)=d^N(x_0)\otimes x_2+x_0\otimes d^N(x_1)=0
 \label{eq3.7}
 \end{equation}
 for $x_0\in C_0$ and $x_1\in C_1$.\\
 
 This tensor product induces the usual tensor product on the underlying graded modules and is essentially unique under this condition.
 
 \subsection{Monoids in $\calc_N$} 
 
 Equipped with the above tensor product the category $\calc_N$ is a monoidal category.\\
 
A monoid of $\calc_N$ is a $N$-complex $\cala$ equipped with a product 
\[
\mu:\cala\otimes\cala\rightarrow \cala
\]
which is associative.\\
 In view of (\ref{eq3.4}) this means that a monoid of $\calc_N$ is a $N$-complex $\cala$ equipped with an associative product $(x,y)\mapsto xy$ which is bilinear, homogeneous of degree 0 and such that
 \begin{equation}
 d(xy)=d(x)y+q^nx d(y)
 \label{eq3.8}
 \end{equation}
 for $x\in \cala^n$ and $y\in \cala$. Thus the monoids of $\calc_N$ are exactly the graded $q$-differential algebras of \cite{mdv:1998a}, (see also in \cite{mdv-ker:1996}).
 
 \section{Concluding remarks}
 
 Given the category $\calc_N$ of $N$-complexes with its monoidal structure described above, the monoids of $\calc_N$ are a natural ``$N$-generalization" of graded differential algebras and coincide with the graded $q$-differential algebras of 
\cite{mdv:1998a}.\\

There are however other natural examples of graded algebras equipped with a $N$-differential such as those considered in Section 7 of \cite{mdv-hen:2002}.\\
In particular, let $N\geq 2$ and $N\geq n\geq 1$ and consider the $N$-homogeneous real algebra $\cala$ generated by $n$ elements $\theta^\lambda$ ($\lambda\in \{1,\dots,n\}$) with relations
\[
\sum_{p\in \cals_N} \theta^{\lambda_{p(1)}}\dots \theta^{\lambda_{p(N)}}=0
\]
for $\lambda_1,\dots,\lambda_N\in \{1,\dots,n\}$ where $\cals_N$ is the group of permutations of $\{1,\dots,N\}$. This algebra is a connected graded algebra $\cala=\oplus_{p\in \mathbb N}\cala^p$ with $\cala^0=\mathbb R\bbbone$ and $\cala^1=\oplus_\lambda\mathbb R\theta^\lambda$. The notion of Koszulity for $N$-homogeneous algebra has been defined in \cite{ber:2001a} where it was shown that $\cala$ is in fact a Koszul algebra as well as its Koszul dual $\cala^!$ 
\cite{ber-mdv-wam:2003}, \cite{ber-mar:2006}. Let us define the graded algebra $\cala(\mathbb R^n)$ by
\[
\cala(\mathbb R^n)=\cala\otimes C^\infty(\mathbb R^n)
\]
where $C^\infty(\mathbb R^n)$ is the algebra of smooth functions on $\mathbb R^n$ and define the homogeneous endomorphism $d$ of degree 1 of $\cala(\mathbb R^n)$ by 
\[
d(a\otimes f)=(-1)^n a\theta^\lambda \otimes \partial_\lambda f
\]
for $a\in \cala^n$ and $f\in C^\infty(\mathbb R^n)$ where the $\partial_\lambda f(x)=\frac{\partial f}{\partial x^\lambda}(x)$ are the partial derivatives with respect to the canonical coordinates $x^\lambda$ of $\mathbb R^n$. One clearly has $d^N=0$ so $d$ is a $N$-differential. Thus $\cala(\mathbb R^n)$ is a graded algebra equipped with a $N$-differential. For $N=2$, $\cala(\mathbb R^n)$ is the graded differential algebra $\Omega(\mathbb R^n)$ of smooth differential forms on $\mathbb R^n$. However for $N\geq 3$, there is no analog of the Leibniz rule on $\cala(\mathbb R^n)$ except for the restriction in degree 0
\[
d:C^\infty(\mathbb R^n)\rightarrow \cala^1(\mathbb R^n)
\]
which is of course a derivation. This comes from the fact that
 $\cala$ has no relation of degree $<N$ and, in particular no relation of degree 2. Neverthesess that kind of generalization of differential forms occurs naturally in lots of problems.\\
 
Other similar graded algebras equipped with $N$-differential (with $N\geq 3$) which are not graded $q$-differential algebras and which are natural are described in \cite{mdv-hen:2002} and in \cite{ber-mdv-wam:2003}.\\

Let us mention a drawback of the notion of graded $q$-differential algebra noticed in \cite{sit:1998} which is that the tensor product of two graded $q$-differential algebras is not a graded $q$-differential algebra but is again simply a graded algebra equipped with a $N$-differential. This is probably connected with the above remarks.\\

These facts which suggest that one needs something more general that the notion of graded $q$-differential algebra are worth noticing in conclusion.

\newpage
\bibliographystyle{plain}

\end{document}